\documentclass[10pt]{article}
   \usepackage{amsmath,amsfonts,amsthm,amssymb,amscd,enumerate, url}
   \usepackage[pdftex]{graphicx,color}  
\theoremstyle{plain}
\newtheorem{thm}{Theorem}[section]
\newtheorem{prop}[thm]{Proposition}
\newtheorem{lemma}[thm]{Lemma}
\newtheorem{cor}[thm]{Corollary}

\theoremstyle{definition}
\newtheorem{defn}[thm]{Definition}

\theoremstyle{remark}
\newtheorem{rmk}[thm]{Remark}

\numberwithin{equation}{section}
\numberwithin{table}{section}
\numberwithin{figure}{section}

\newcommand{\vol}{\ensuremath{\mathrm{vol}}}
\newcommand{\real}{\ensuremath{\mathrm{Re}}}
\newcommand{\imag}{\ensuremath{\mathrm{Im}}}
\newcommand{\spn}{\ensuremath{\mathrm{span}}}

\newcommand{\G}{\ensuremath{\mathrm{G}_2}}
\newcommand{\SP}{\ensuremath{\mathrm{Spin}(7)}}

\newcommand{\SUth}{\ensuremath{\mathrm{SU}(3)}}

\newcommand{\ph}{\ensuremath{\varphi}}
\newcommand{\ps}{\ensuremath{\psi}}
\newcommand{\Ph}{\ensuremath{\Phi}}

\newcommand{\R}{\ensuremath{\mathbb R}}
\newcommand{\C}{\ensuremath{\mathbb C}}
\newcommand{\Q}{\ensuremath{\mathbb H}}

\newcommand{\Oc}{\ensuremath{\mathbb O}}

\newcommand{\st}{\ensuremath{\ast}}
\newcommand{\hk}{\mathbin{\! \hbox{\vrule height0.3pt width5pt depth 0.2pt \vrule height5pt width0.4pt depth 0.2pt}}}

\newcommand{\oi}{\ensuremath{\mathbf i}}
\newcommand{\oj}{\ensuremath{\mathbf j}}
\newcommand{\ok}{\ensuremath{\mathbf k}}
\newcommand{\oee}{\ensuremath{\mathbf e}}

\newcommand{\bOm}{\ensuremath{\overline{\Omega}}}
\newcommand{\curl}{\ensuremath{\operatorname{curl}}}
\newcommand{\grad}{\ensuremath{\operatorname{grad}}}
\newcommand{\dive}{\ensuremath{\operatorname{div}}}
\newcommand{\nab}[1]{\ensuremath{\nabla_{\! \! #1 \,}}}
\newcommand{\dirac}{\ensuremath{\not \! \! \mathcal D}}
\newcommand{\cross}{\ensuremath{\! \times \!}}

\begin{document}

\title{Some Notes on $\G$ and $\SP$ Geometry}

\author{Spiro Karigiannis \\ Department of Pure Mathematics, University of Waterloo \\ \tt{karigiannis@math.uwaterloo.ca}}

\date{\today}

\maketitle

\begin{abstract}
We collect together various facts about $\G$ and $\SP$ geometry which are likely well known but which do not seem to have appeared explicitly in the literature before. These notes should be useful to graduate students and new researchers in $\G$ and $\SP$ geometry.
\end{abstract}

\section{Introduction} \label{introsec}

These notes consist of a collection of facts about manifolds with $\G$ or $\SP$ structures.

\begin{itemize}
\item In Section~\ref{signssec} we discuss the various sign and orientation conventions commonly used in $\G$ and $\SP$ geometry and the relations between them.
\item In Section~\ref{formssec} we discuss the relations between the form type decompositions of Calabi-Yau $3$-folds and $\G$ manifolds. If $X^6$ is a Calabi-Yau $3$-fold, then $M^7 = S^1 \times X^6$, with the product metric, is a $\G$ manifold, and we can decompose forms on $M^7$ into components determined by $\G$ representation theory. These can be compared to the
$\SUth$ type decompositions arising from the Calabi-Yau geometry. We find that there is a non-trivial mixing between the K\"ahler moduli and complex moduli of $X^6$ inside the $\G$ moduli of $S^1 \times X^6$.
\item In Section~\ref{operatorssec} we discuss some natural first order differential operators on manifolds with $\G$ structure: the gradient, divergence, curl, and Dirac operator. We also discuss the relations between these operators.
\end{itemize}

The reader is expected to be familiar with $\G$ and $\SP$ structures. Some references for $\G$ and $\SP$ structures are~\cite{Br1, Br2, BS, HL, J1, J2, K1}.

{\em Disclaimer.} These notes are purely expository and are made available as a potentially useful resource to others. In particular the reference list at the end is incomplete and any particular inclusions and/or omissions are not in any way meant to imply anything.

\section{Signs and orientations in $\G$ and $\SP$ Geometry} \label{signssec}  

\subsection{$\G$ structures} \label{g2sec}

The differential forms that describe $\G$ and $\SP$ structures can be defined in terms of octonion algebra. If we define the multiplication on the octonions $\Oc = \Q \oplus \Q \oee = \R^8$ via the
Cayley-Dickson process~\cite{HL}, we have
\begin{equation*}
(a + b \oee) \cdot (c + d \oee) = (a c - \bar d b) + (d a + b \bar c)
\oee \qquad a, b, c, d \in \Q
\end{equation*}
in terms of quaternion multiplication. Let $\langle \cdot, \cdot \rangle$ denote the standard Euclidian inner product on $\R^8$. Following~\cite{HL}, on $\imag (\Oc) = \R^7$ we define the $3$-form $\ph$ by
\begin{equation*}
\ph(x,y,z) \, = \, \langle x, y z \rangle \qquad x,y,z \in \imag (\Oc),
\end{equation*}
and its dual $4$-form $\ps$ by
\begin{equation*}
\ps(x,y,z,w) \, = \, \frac{1}{2} \langle x, [y, z, w] \rangle \qquad x,y,z,w \in \imag (\Oc),
\end{equation*}
where $[x,y,z] = (xy)z - x(yz)$ is the associator.

In terms of the standard basis for $\R^8 = \Oc$ we have the coordinates $x^0, x^1, x^2, x^3$, $y^0, y^1, y^2, y^3$ where the $x^i$'s are coordinates on $\Q$ and the $y^j$'s are coordinates on $\Q \oee$. We
take the orientation given by the volume form
\begin{equation*}
\vol_8 \, = \, dx^0 dx^1 dx^2 dx^3 dy^0 dy^1 dy^2 dy^3.
\end{equation*}
Here and in the sequel, we omit the wedge broduct symbol $\wedge$ to avoid notational clutter. The induced orientation on $\R^7 = \imag (\Oc)$ is taken to be the one given by
\begin{equation*}
\vol_7 \, = \, dx^1 dx^2 dx^3 dy^0 dy^1 dy^2 dy^3.
\end{equation*}
With respect to this orientation, the $4$-form $\ps$ is the Hodge dual (on $\R^7$) of $\ph$:
\begin{equation*}
\psi = \st_7 \ph.
\end{equation*}

In these coordinates, the forms $\ph$ and $\ps$ can be written as
\begin{equation} \label{convention1pheq}
\begin{aligned}
\ph \, =  \, \, & dx^1 dx^2 dx^3 - dx^1 dy^2 dy^3 - dy^1 dx^2 dy^3 - dy^1 dy^2
dx^3 \\ & {} - dy^0 dx^1 dy^1 - dy^0 dx^2 dy^2 - dy^0 dx^3 dy^3,
\end{aligned}
\end{equation}
and
\begin{equation} \label{convention1pseq}
\begin{aligned}
\ps \, = \, \, & dy^0 dy^1 dy^2 dy^3 - dy^0 dy^1 dx^2 dx^3 - dy^0 dx^1 dy^2 dx^3
- dy^0 dx^1 dx^2 dy^3 \\ & {} - dx^2 dy^2 dx^3 dy^3 - dx^3 dy^3 dx^1
dy^1 - dx^1 dy^1 dx^2 dy^2.
\end{aligned}
\end{equation}
It is easy to see that $\ph \wedge \ps = 7 \vol_7$.

There are two useful ways to interpret and remember these expressions. First, we can consider $\R^7 = \C^3 \oplus \R$ where we take $z^j = x^j + i y^j$ as complex coordinates for $j = 1,2,3$ and $y^0 = t$ is a real coordinate. Then the standard orientation on $\C^3$ is
\begin{equation*}
\vol_6 \, = \, dx^1 dy^1 dx^2 dy^2 dx^3 dy^3,
\end{equation*}
and we have
\begin{equation*}
\vol_7 \, = \, dt \wedge \vol_6.
\end{equation*}
With this identification, the forms become
\begin{align*}
\ph \, = \, \, & \real (\Omega_3) - dt \wedge \omega_3, \\ \ps \, = \, \, & - dt \wedge \imag (\Omega_3) - \frac{1}{2}\omega_3^2,
\end{align*}
where
\begin{align*}
\Omega_3 \, = \, \, & dz^1 dz^2 dz^3, \\ \omega_3 \, = \, \, & \frac{i}{2} (dz^1 d\bar{z}^1 + dz^2 d\bar{z}^2 + dz^3 d\bar{z}^3) = dx^1 dy^1 + dx^2 dy^2 + dx^3 dy^3,
\end{align*}
are the standard holomorphic volume form and K\"ahler form on $\C^3$, respectively.

A second way to interpret and remember these forms is to consider $\R^7 = \imag (\Oc) = \imag (\Q) \oplus \Q \oee = \R^3 \oplus \R^4$ where we take the $x^i$'s as coordinates on $\R^3$ for $i = 1,2,3$ and the $y^j$'s as coordinates on $\R^4$ for $j = 0,1,2,3$. We take the orientations
\begin{equation*}
\vol_3 = dx^1 dx^2 dx^3 \qquad \qquad \text{and} \qquad \qquad \vol_4 = dy^0 dy^1 dy^2 dy^3
\end{equation*}
on $\R^3$ and $\R^4$, respectively. This gives us the following basis of \emph{anti-self dual} $2$-forms on $\R^4$:
\begin{align*}
\eta^-_1 \, = \, \, & dy^0 dy^1 - dy^2 dy^3, \\ \eta^-_2 \, = \, \, & dy^0 dy^2 - dy^3 dy^1, \\ \eta^-_3 \, = \, \, & dy^0 dy^3 - dy^1 dy^2.
\end{align*}
With this identification, the forms are
\begin{align*}
\ph \, = \, \, & \vol_3 + \, dx^1 \eta^-_1 + dx^2 \eta^-_2 + dx^3 \eta^-_3, \\ \ps \, = \, \, & \vol_4 - dx^1 dx^2 \eta^-_3 - dx^2 dx^3 \eta^-_1 - dx^3 dx^1 \eta^-_2.
\end{align*}
We can identify $\R^7 = \Lambda^2_- (\R^4) = \R^3 \oplus \R^4$ and in this way the canonical form $\theta$ on the total space is given by
\begin{equation*}
\theta \, = \, x^1 \eta^-_1 + x^2 \eta^-_2 + x^3 \eta^-_3,
\end{equation*}
and hence
\begin{equation*}
\ph = \vol_3 + d\theta.
\end{equation*}

To extract the metric from the $3$-form $\ph$, we use the identity
\begin{equation*}
(u \hk \ph) \wedge (v \hk \ph) \wedge \ph \, = \, -6 \langle u, v \rangle \vol_7
\end{equation*}
The forms $\ph$ and $\ps$ are characterized by the fact that they are invariant under the action of the group $\G \subset SO(7)$.

Alternatively, we can instead consider the opposite orientation $\widetilde\vol_7 = - \vol_7$. This is equivalent to keeping the orientation fixed and changing the sign of one coordinate. (And changing the sign of the Hodge star $\st$ operator.) Let us change $y^0 \to - y^0$. Then in terms of $\C^3 \oplus \R$, the forms become
\begin{align*}
\ph \, = \, \, & \real (\Omega_3) + dt \wedge \omega_3, \\ \ps \, = \, \, & - dt \wedge \imag (\Omega_3) + \frac{1}{2} \omega_3^2.
\end{align*}
Using this convention, in terms of $\R^3 \oplus \R^4$, the forms are
\begin{align*}
\ph \, = \, \, & \vol_3 - \, dx^1 \eta^+_1 - dx^2 \eta^+_2 - dx^3 \eta^+_3, \\ \ps \, = \, \, & \vol_4 - \, dx^1 dx^2 \eta^+_3 - dx^2 dx^3 \eta^+_1 - dx^3 dx^1 \eta^+_2,
\end{align*}
where
\begin{align*}
\eta^+_1 \, = \, \, & dy^0 dy^1 + dy^2 dy^3, \\ \eta^+_2 \, = \, \, & dy^0 dy^2 + dy^3 dy^1, \\ \eta^+_3 \, = \, \, & dy^0 dy^3 + dy^1 dy^2,
\end{align*}
is the corresponding basis of \emph{self-dual} $2$-forms on $\R^4$.

In this convention, the metric is obtained from
\begin{equation*}
(u \hk \ph) \wedge (v \hk \ph) \wedge \ph \, = \, +6 \langle u, v \rangle \vol_7
\end{equation*}
In this case the model space for the $\G$ structure is $\Lambda^2_+ (\mathbb R^4)$.

\begin{rmk}
If instead we defined the $3$-form $\tilde \ph$ by
\begin{equation*}
\tilde \ph \, = \, \vol_3 - \sum_k dx^k \eta^-_k \qquad \qquad \text{or} \qquad \qquad \tilde \ph \, = \, \vol_3 + \sum_k dx^k \eta^+_k
\end{equation*}
then these forms would determine a metric of signature $(3,4)$ or $(4,3)$ and correspond to the non-compact split form of $\G$.
\end{rmk}

\subsection{$\SP$ structures} \label{spin7sec}

We now move on to \SP\ structures. Following~\cite{HL}, there is a $3$-fold cross product on $\R^8 = \Oc$ defined via octonion multiplication by
\begin{equation*}
X(x,y,z) \, = \, \frac{1}{2} (x (\bar{y} z) - z ( \bar{y} x) ),
\end{equation*}
which allows us to define the canonical $4$-form $\Ph$ on $\Oc$ as
\begin{equation*}
\Ph(x,y,z,w) \, = \, \langle x, X(y,z,w) \rangle.
\end{equation*}
In terms of the coordinates $x^i, y^i$, $i=0,1,2,3$, the form $\Ph$ can be written as
\begin{align*}
\Ph \, = \, & dx^0 dx^1 dx^2 dx^3 - dx^0 dx^1 dy^2 dy^3 - dx^0 dy^1 dx^2 dy^3 - dx^0 dy^1 dy^2 dx^3 \\ \, & {} - dx^0 dy^0 dx^1 dy^1 - dx^0 dy^0 dx^2 dy^2 - dx^0 dy^0 dx^3 dy^3 \\ \, & {} + dy^0 dy^1 dy^2 dy^3 - dy^0 dy^1 dx^2 dx^3 - dy^0 dx^1 dy^2 dx^3 - dy^0 dx^1 dx^2 dy^3 \\ \, & {} - dx^2 dy^2 dx^3 dy^3 - dx^3 dy^3 dx^1 dy^1 - dx^1 dy^1 dx^2 dy^2
\end{align*}
which can be more succintly written as
\begin{equation*}
\Ph \, = \, dx^0 \wedge \ph + \psi.
\end{equation*}
With respect to the standard orientation $\vol_8$, the $4$-form $\Ph$ is self-dual: $\st \Ph = \Ph$, and it satisfies $\Ph^2 = 14 \vol_8$. The above expression corresponds to the model space $\Oc =
\R^8 = \R \oplus \R^7$.

There are two more ways to interpret and remember the $4$-form $\Ph$. First, one can identify $\Oc = \R^8 = \C^4$ and take complex coordinates $z^j = x^j + i y^j$ for $j = 0,1,2,3$. Note that our
chosen orientation $\vol_8$ is equal to the canonical orientation determined by the complex structure, specifically
\begin{equation*}
\vol_8 \, = \, dx^0 dy^0 dx^1 dy^1 dx^2 dy^2 dx^3 dy^3.
\end{equation*}
In these coordinates, the $4$-form $\Ph$ can be written as
\begin{equation*}
\Ph \, = \, \real (\Omega_4) - \frac{1}{2}\omega_4^2,
\end{equation*}
where
\begin{align*}
\Omega_4 \, = \, \, & dz^0 dz^1 dz^2 dz^3, \\ \omega_4 \, = \, \, & \frac{i}{2} (dz^0 d\bar{z}^0 + dz^1 d\bar{z}^1 + dz^2 d\bar{z}^2 + dz^3 d\bar{z}^3) \\ = \, \, & dx^0 dy^0 + dx^1 dy^1 + dx^2 dy^2 + dx^3 dy^3,
\end{align*}
are the standard holomorphic volume form and K\"ahler form on $\C^4$, respectively.

Alternatively, if we use the description of the $\G$ forms $\ph$ and $\ps$ in terms of the decomposition $\R^7 = \R^3 \oplus \R^4$, then we get a description of $\Ph$ in terms of $\R^8 = \R^4 \oplus \R^4$. Let us define $\vol_X = dx^0 dx^1 dx^2 dx^3$ and $\vol_Y = dy^0 dy^1 dy^2 dy^3$. Let $\eta^-_i = dy^0 dy^i - dy^j dy^k$ be a basis of anti-self-dual $2$-forms on the copy of $\R^4$ with coordinates $y^i$, where $i,j,k$ is a cyclic permutation of $1,2,3$. Similarly define $\beta^-_i = dx^0 dx^i - dx^j dx^k$ to be anti-self dual $2$-forms on the other copy of $\R^4$. Then we have
\begin{equation*}
\Ph \, = \,  \vol_X + \beta^-_1 \eta^-_1 + \beta^-_2 \eta^-_2 + \beta^-_3 \eta^-_3 + \vol_Y.
\end{equation*}
Before interpreting this decomposition $\R^8 = \R^4 \oplus \R^4$ in terms of spinor spaces, we remark that if we chose the other convention for $\ph$ and $\ps$, which corresponds to a change or orientation on $\R^7$ (and hence on $\R^8$), we would obtain the following two forms for $\Ph$:
\begin{align*}
\Ph \, = \, \, & \real (\Omega_4) + \frac{1}{2}\omega_4^2, \\ \Ph \, = \, \, & \vol_X - \beta^+_1 \eta^+_1 - \beta^+_2 \eta^+_2 - \beta^+_3 \eta^+_3 + \vol_Y,
\end{align*}
where the $\eta^+_i$ and $\beta^+_i$ are now self-dual $2$-forms on the two copies of $\R^4$. The $4$-form $\Ph$ is still self-dual.

Consider now the space $\mathbb S_- (\R^4)$ of negative chirality spinors over $\R^4$. This is a quaternionic line bundle over $\R^4$. Let $e^0, e^1, e^2, e^3$ be an oriented orthonormal basis of $\R^4$. Define $\omega^-_i = e^0 \wedge e^i - e^j \wedge e^k$ where $i,j,k$ is a cyclic permutation of $1,2,3$. It is easy to check that under Clifford multiplication, $\omega^-_i \cdot \omega^-_i = -2(1 + \vol)$ and $\omega^-_i \cdot \omega^-_j = - 2 \omega^-_k$ where $\vol = e^0 \cdot e^1 \cdot e^2 \cdot e^3$ is the volume form. Since for dimension $n=4$, Clifford multiplication by $\gamma = -\vol$ is equal to $\pm 1$ on the spinor spaces $\mathbb S_{\pm}$, we see that
\begin{equation*}
\frac{1}{2}\omega^-_1 , \quad \quad \frac{1}{2}\omega^-_2, \quad \quad \frac{1}{2}\omega^-_3,
\end{equation*}
act as right multiplication by the quaternions $\oi$, $\oj$, $\ok$, respectively, on the fibre $\mathbb S_- = \Q$. Thus, for any choice of unit spinor $s_0 \in \mathbb S_+$, we obtain an orthonormal basis
\begin{equation*}
s_0, \qquad s_1 = \frac{1}{2}\omega^-_1 \cdot s_0, \qquad s_2 = \frac{1}{2}\omega^-_2 \cdot s_0, \qquad s_3 = \frac{1}{2}\omega^-_3 \cdot s_0,
\end{equation*}
of $\mathbb S_-$ and conversely every orthonormal basis can be written in this way. Following~\cite{DH}, we can multiply two spinors $s_i$, $s_j$ to obtain an endomorphism $s_i \circ s_j$ of $\mathbb S_-$, which is itself Clifford multiplication by some form. That is, the product of two spinors is a form. The product is defined by
\begin{equation*}
(s_i \circ s_j) s_k \, = \, \langle s_j, s_k \rangle s_i.
\end{equation*}
Then it is easy to check explicitly using an orthonormal basis $s_0$, $s_1$, $s_2$, $s_3$, that as endomorphisms of $\mathbb S_-$, we have
\begin{equation*}
\frac{1}{2}\omega^-_i \, = \, - (s_0 \circ s_i - s_j \circ s_k) + (s_i \circ s_0 - s_k \circ s_j )
\end{equation*}
where $i,j,k$ is a cyclic permutation of $1,2,3$. Now suppose we choose some different orthonormal basis $\tilde s_k = A^l_k s_l$, where $k, l = 0,1,2,3$. Then the associated orthonormal basis $\{
\beta^-_i\}$ of anti-self-dual $2$-forms on the fibre $\mathbb S_-$ will change by the linear map $\Lambda^2_- (A^*)$, where $A^*$ is the inverse transpose of $A$ and $\Lambda^2_- (A^*)$ means the induced linear operator on $\Lambda^2_- (\mathbb S_-)^*$. But by the above expression, the associated orthonormal basis $\{\eta^-_i\}$ of anti-self dual $2$-forms on the base $\R^4$ will change by $\Lambda^2_- (A)$. Thus the expression
\begin{equation*}
\beta^-_1 \eta^-_1 + \beta^-_2 \eta^-_2 + \beta^-_3 \eta^-_3
\end{equation*}
is independent of choice of orthonormal basis and is hence a well-defined $4$-form on the total space $\mathbb S_- (\R^4)$. Therefore with this choice of convention, the model space for the $\SP$ structure on $\R^8$ is $\mathbb S_- (\R^4)$.

If instead we had used the other sign/orientation convention, we would be using self-dual $2$-forms on the base and fibre, and the total space would be the space of positive chirality spinors, $\mathbb S_+ (\R^4)$.

\section{Relationship between $\G$ manifolds and Calabi-Yau $3$-folds} \label{formssec}

\subsection{$\G$ structures on $S^1 \times X^6$} \label{g2cy3}

Let $X^6$ be a Calabi-Yau $3$-fold with K\"ahler form $\omega$ and non-vanishing holomorphic $(3,0)$ form $\Omega$. The following relations hold:
\begin{align*}
& \omega \wedge \Omega \, = \, \omega \wedge \bOm \, = \, \omega \wedge \real(\Omega) \, = \, \omega \wedge \imag(\Omega) \, =\,  0, \\ & \frac{\omega^3}{6} \, = \, \vol_6 \, = \, \frac{i}{8} \Omega \wedge \bOm,\qquad \qquad {|\real(\Omega)|}^2 \, = \, {|\imag(\Omega)|}^2 \, = \, 4, \\ & \real(\Omega) \wedge \imag(\Omega) \, = \, \frac{1}{4i}(\Omega + \bOm) \wedge (\Omega - \bOm) \, = \, \frac{i}{2} \Omega \wedge \bOm \, = \, 4 \vol_6, \\ & \st_6 \real(\Omega) \, = \, \imag(\Omega), \qquad \qquad \st_6 \imag(\Omega) \, = \, - \real(\Omega),
\end{align*}
where $\vol_6$ and $\st_6$ are the volume form and Hodge star operator on $X^6$, respectively. Also, $| \cdot |^2$ is the pointwise norm on forms on $X^6$.

Now let $t$ be an angle coordinate for the circle $S^1$, so $dt$ is the globally defined volume form on $S^1$ with respect to the standard round metric. With the product metric on $M^7 = S^1 \times X^6$, the
$7$-manifold has holonomy contained in \G. (In fact the holonomy is \SUth.) We can take the associated $\G$ $3$-form $\ph$ to be
\begin{equation*}
\ph \, = \, \real(\Omega) - dt \wedge \omega,
\end{equation*}
which determines the dual $4$-form $\ps = \st_7 \ph$ as
\begin{equation*}
\ps \, = \, -dt \wedge \imag(\Omega) - \frac{\omega^2}{2},
\end{equation*}
where we use $\st_7$ to denote the Hodge star operator on $M^7$ and $\vol_7$ the volume form on $M^7$. We see that
\begin{align*}
\ph \wedge \ps \, = \, \, & dt \wedge \real(\Omega) \wedge \imag(\Omega) + dt \wedge \frac{\omega^3}{2} \\ = \, \, & 4 dt \wedge \vol_6 + 3 dt \wedge \vol_6 \\ = \, \, & 7 dt \wedge \vol_6 = 7 \vol_7 = {|\ph|}^2 \vol_7
\end{align*}
as expected.

\subsection{Form type decompositions} \label{formstypesec}

On the Calabi-Yau $3$-fold $X^6$ we can decompose the complex valued differential forms into $(p,q)$ types given by the complex structure. More specifically, the complex valued $2$-forms and $3$-forms
decompose as
\begin{align*}
\Omega^2(X^6, \C) \, = \, \, & \Omega^{2,0} \oplus \Omega^{1,1} \oplus \Omega^{0,2}, \\ \Omega^3(X^6, \C) \, = \, \, & \Omega^{3,0} \oplus \Omega^{2,1} \oplus \Omega^{1,2} \oplus \Omega^{0,3},
\end{align*}
where $\Omega^{p,q}$ are the complex-valued forms of type $(p,q)$ and $\Omega^{q,p} = \overline{\Omega^{p,q}}$. Since $X^6$ is K\"ahler, the $(1,1)$ forms further decompose into
\begin{equation*}
\Omega^{1,1} \, = \, \spn(\omega) \oplus \Omega^{1,1}_0
\end{equation*}
where $\Omega^{1,1}_0$ are the $(1,1)$ forms which are pointwise orthogonal to the K\"ahler form $\omega$, and $\spn(\omega) = \{ f \omega; f \in C^{\infty}(X^6) \}$.

On the $\G$ manifold $M^7$, the complex valued $2$-forms and $3$-forms decompose as
\begin{align*}
\Omega^2 (M^7, \C) \, = \, \, & \Omega^2_7 \oplus \Omega^2_{14}, \\ \Omega^3 (M^7, \C) \, = \, \, & \Omega^3_1 \oplus \Omega^3_7 \oplus \Omega^3_{27},
\end{align*}
where the subspaces $\Omega^k_l$ are defined by
\begin{align*}
\Omega^2_7 \, = \, \, & \{ X \hk \ph; X \in \Gamma(TM^7) \} = \{ \beta \in \Omega^2; \st_7 (\ph \wedge \beta) = -2 \beta \}, \\ \Omega^2_{14} \, = \, \, & \{ \beta \in \Omega^2; \beta \wedge \ps = 0 \} = \{ \beta \in \Omega^2; \st_7 (\ph \wedge \beta) = \beta \}, \\ \Omega^3_1 \, = \, \, & \{ f \ph ; f \in C^{\infty}(M^7) \}, \\ \Omega^3_7 \, = \, \, & \{ X \hk \ps; X \in \Gamma(TM^7) \}, \\ \Omega^3_{27} \, = \, \, & \{ \eta \in \Omega^3; \eta \wedge \ph = 0 \, \text{ and } \, \eta \wedge \ps = 0 \}.
\end{align*}

\begin{rmk} \label{signsrmk}
If we had instead used the opposite orientation, we would have $\ph = \real(\Omega) + dt \wedge \omega$, with dual $4$-form $\ps = - dt \wedge \imag(\Omega) + \frac{\omega^2}{2}$. With this convention, the eigenspaces $\Omega^2_7$ and $\Omega^2_{14}$ of the operator $\beta \mapsto \st_7 (\ph \wedge \beta)$ on $\Omega^2$ would correspond to eigenvalues $+2$ and $-1$, respectively.
\end{rmk}

\subsubsection*{The case of $3$-forms}

Let $\eta$ be a complex valued $3$-form on $M^7$. We can decompose it as
\begin{equation*}
\eta \, = \, \eta_3 + dt \wedge \eta_2,
\end{equation*}
where $\eta_2$ and $\eta_3$ are a $2$-form and a $3$-form on $X^6$, respectively. Note that strictly speaking, these forms also depend on the parameter $t$. This will not affect our computations, however. We simply treat forms on $X^6$ as depending on a smooth angular parameter $t$.

We write $\eta_3 = \eta_{3,0} + \eta_{2,1} + \eta_{1,2} + \eta_{0,3}$ and $\eta_2 = \eta_{2,0} + \eta^0_{1,1} + f \omega + \eta_{0,2}$, decomposing them into types determined by the K\"ahler structure of $X^6$.

\begin{prop} \label{3-27-prop}
The $3$-form $\eta$ is in $\Omega^3_{27}$ if and only if the following equations are satisfied:
\begin{align*}
& \omega \wedge \eta_{2,1} + \frac{1}{2} \Omega \wedge \eta_{0,2} \, = \, 0, \\ & \omega \wedge \eta_{1,2} + \frac{1}{2} \bOm \wedge \eta_{2,0} \, = \, 0, \\ & \frac{1}{2} \bOm \wedge \eta_{3,0} + \frac{1}{2} \Omega
\wedge \eta_{0,3} \, = \, 0, \\ & \frac{i}{2} \bOm \wedge \eta_{3,0} - \frac{i}{2} \Omega \wedge \eta_{0,3} + \frac{f}{2} \omega^3 \, = \, 0.
\end{align*}
\end{prop}
\begin{proof}
First, we compute $\eta \wedge \ph$ and decompose into types:
\begin{align*}
\eta \wedge \ph \, = \, \, & \left( \eta_{3,0} + \eta_{2,1} + \eta_{1,2} + \eta_{0,3} + dt \wedge (\eta_{2,0} + \eta^0_{1,1} + f \omega + \eta_{0,2}) \right) \wedge \left( \frac{1}{2} ( \Omega + \bOm ) - dt \wedge \omega \right) \\ = \, \, & \frac{1}{2} \eta_{3,0} \wedge \bOm + \frac{1}{2} \eta_{0,3} \wedge \Omega + dt \wedge \omega \wedge \eta_{2,1} + dt \wedge \omega \wedge \eta_{1,2} + \frac{1}{2} dt \wedge \eta_{2,0} \wedge \bOm + \frac{1}{2} dt \wedge \eta_{0,2} \wedge \Omega,
\end{align*}
where all other terms are zero due to type considerations. Now collecting terms of the same type: $(3,3)$, $dt \wedge (3,2)$, and $dt \wedge (2,3)$, gives the first three equations above. Similarly we
compute
\begin{align*}
\eta \wedge \ps \, = \, \, & \left( \eta_{3,0} + \eta_{2,1} + \eta_{1,2} + \eta_{0,3} + dt \wedge (\eta_{2,0} + \eta^0_{1,1} + f \omega + \eta_{0,2}) \right) \wedge \left( -dt \wedge \frac{1}{2i} ( \Omega - \bOm ) - \frac{\omega^2}{2} \right) \\ = \, \, & \frac{i}{2} dt \wedge \eta_{3,0} \wedge \bOm - \frac{i}{2} dt \wedge \eta_{0,3} \wedge \Omega - \frac{f}{2} dt \wedge \omega^3,
\end{align*}
where again all other terms are zero due to type considerations, and we have also used the fact that $\eta^0_{1,1} \wedge \omega^2$ is zero, which follows from the fact that $\eta^0_{1,1} \perp \omega$ and
$\st_6 \omega = \frac{\omega^2}{2}$. This expression is all of type $dt \wedge (3,3)$, and setting it equal to zero gives the fourth equation above.
\end{proof}

As an example, if we take $\eta = \real(\Omega) + f dt \wedge \omega$ for some function $f$, we can check easily that this $\eta$ is in $\Omega^3_{27}$ if and only if $f = \frac{4}{3}$. That is,
\begin{equation*}
\real(\Omega) + \frac{4}{3} dt \wedge \omega \, \,  \in \, \,  \Omega^3_{27}.
\end{equation*}
Recall that by definition we have
\begin{equation*}
\real(\Omega) - dt \wedge \omega \, \, \in \, \, \Omega^3_1,
\end{equation*}
from which it follows immediately that the real $3$-forms $\real(\Omega)$ and $dt \wedge \omega$ are both in $\Omega^3_1 \oplus \Omega^3_{27}$, and not in a strictly smaller subspace. We will see shortly that $\imag(\Omega)$ lies in $\Omega^3_7$.

Before we move on to $\Omega^3_7$, consider the case when $\eta$ is real.  Then $\overline{\eta} = \eta$, so it follows that $\eta_{q,p} = \overline{\eta_{p,q}}$. Now suppose that we have $\eta = \eta_{3,0} + \overline{\eta_{3,0}} + \eta_{2,1} + \overline{\eta_{2,1}} + dt \wedge \left( \eta_{2,0} + \overline{\eta_{2,0}} + f \omega + \eta^0_{1,1} \right)$ is some real $3$-form, where $f$ and $\eta^0_{1,1}$ are both
real. Necessarily $\eta_{3,0} = g \Omega$ for some function $g$. Substituting this expression into the equations of Proposition~\ref{3-27-prop}, we find easily that $g = \frac{8}{3} f$ and is thus also real. Therefore we have
\begin{cor} \label{3-27-real-cor}
The {\em real} $3$-forms of type $\Omega^3_{27}$ are given by:
\begin{itemize} \setlength\itemsep{-1mm}
\item $\spn ( \real(\Omega) + \frac{4}{3} dt \wedge \omega )$, which is (pointwise) $1$-dimensional.
\item all real forms of type $dt \wedge (1,1)_0$ (where $(1,1)_0$ are the $(1,1)$ forms orthogonal to $\omega$), which is (pointwise) $8$-dimensional.
\item the real $3$-forms $\eta_{2,1} + \overline{\eta_{2,1}} + dt \wedge ( \eta_{2,0} + \overline{\eta_{2,0}} )$, where $\frac{1}{2} \Omega \wedge \overline{\eta_{2,0}} + \omega \wedge \eta_{2,1} = 0$,
which is (pointwise) $18$-dimensional. This is because for any real $(2,1)$ form $\eta_{2,1} + \overline{\eta_{2,1}}$, the equation $\frac{1}{2} \Omega \wedge \overline{\eta_{2,0}} + \omega \wedge \eta_{2,1} = 0$ can be solved uniquely for $\eta_{2,0}$ due to the non-degeneracy of the holomorphic volume form $\Omega$.
\end{itemize}
\end{cor}
This gives the expected result for the (pointwise) dimension of the space of {\em real} $\Omega^3_{27}$ forms to be $1 + 8 + 18 = 27$.

We now move on to the $3$-forms $\eta$ of type $\Omega^3_7$. These are still orthogonal to $\ph$, so we still require the condition $\eta \wedge \ps = 0$, but now the map $\eta \mapsto \ph \wedge \eta$ is an isomorphism of $\Omega^3_7$ onto $\Omega^6_7$.

\begin{prop} \label{3-7-prop}
The {\em real} $3$-forms of type $\Omega^3_7$ are the following:
\begin{itemize} \setlength\itemsep{-1mm}
\item $\spn( \imag(\Omega) ) \ $, which is (pointwise) $1$-dimensional.
\item the real $3$-forms $\eta_{2,1} + \overline{\eta_{2,1}} + dt \wedge ( \eta_{2,0} + \overline{\eta_{2,0}} )$, where $(\eta_{2,1}, \eta_{2,0})$ is in the orthogonal complement of the kernel of the
linear map
\begin{equation*}
L : (\eta_{2,1}, \eta_{2,0}) \mapsto \frac{1}{2}
\Omega \wedge \overline{\eta_{2,0}} + \omega \wedge \eta_{2,1}.
\end{equation*}
This space is (pointwise) $6$-dimensional.
\end{itemize}
This gives the expected result of $1 + 6 = 7$ for the (pointwise) dimension of the space of real $\Omega^3_7$ forms.
\end{prop}
\begin{proof}
Suppose that the real $3$-form $\eta$ is in $\Omega^3_7$. Since we have already shown that $dt \wedge \omega$ is in $\Omega^3_1 \oplus \Omega^3_{27}$ and that all the forms of the type $dt \wedge (1,1)_0$ are in $\Omega^3_{27}$, we can write that
\begin{equation*}
\eta = \eta_{3,0} + \overline{\eta_{3,0}} + \eta_{2,1} +
\overline{\eta_{2,1}} + dt \wedge ( \eta_{2,0} + \overline {\eta_{2,0}} ).
\end{equation*}
We must have $\eta_{3,0} = g \Omega$ for some function $g$, and substituting this into the fourth equation of Proposition~\ref{3-27-prop}, which holds here since $\Omega^3_7$ forms also satisfy $\eta \wedge \ps = 0$, we get that $g$ must be purely imaginary. This gives
\begin{equation*}
\eta = h \imag(\Omega) + \eta_{2,1} + \overline{\eta_{2,1}} + dt \wedge ( \eta_{2,0} + \overline {\eta_{2,0}} )
\end{equation*}
for some real function $h$. It remains to deal with the other terms. Consider the linear map $L$ from the underlying real vector space of $\Omega^{2,1} \oplus \Omega^{2,0}$ to the underlying real vector space of $\Omega^{3,2}$, given by $L(\eta_{2,1}, \eta_{2,0}) = \frac{1}{2} \Omega \wedge \overline{\eta_{2,0}} + \omega \wedge \eta_{2,1}$. We have shown in Corollory~\ref{3-27-real-cor} that $L$ has an $18$-dimensional kernel, so by the rank-nullity theorem, since the domain is $24$-dimensional, the orthogonal complement to the kernel is mapped isomorphically onto the image, a $6$-dimensional real
vector space. Note that the map $L$ is precisely wedge product with the $\G$ $3$-form $\ph$. Since the $\Omega^3_7$ forms are precisely those which are mapped isomorphically onto $\Omega^6_7$ by wedge product with $\ph$, this gives those $\Omega^3_7$ forms which are mapped to $6$-forms of the form $dt \wedge \alpha_{3,2} + dt \wedge \overline{\alpha_{3,2}}$, a $6$-dimensional space, and the remaining $\Omega^3_7$ form is $\imag(\Omega)$, which is sent to $\real(\Omega) \wedge \imag(\Omega) = 4 \vol_6$ under wedge product with $\ph$.
\end{proof}

Here is another way to describe the real $\Omega^3_7$ forms. Recall that they are given by $X \hk \ps$, where $X$ is a real vector field on $M^7$. We can write this as
\begin{equation*}
X = a^i \frac{\partial}{\partial z^i} + \bar a^i \frac{\partial}{\partial \bar z^i} + h \frac{\partial}{\partial t}
\end{equation*}
in terms of local complex coordinates $z^i$ on $X^6$. Then it is easy to check that
\begin{align*}
& X \hk \ps \, = \, -h \imag(\Omega) + a^i \eta_i + \bar a^i \overline{\eta}_i, \\ \, \, & \text{ where } \ \eta _ i = dt \wedge \frac{1}{2i} \left( \frac{\partial}{\partial z^i} \hk \Omega \right) - \frac{\partial}{\partial z^i} \hk \frac{\omega^2}{2}.
\end{align*}
This gives a canonical basis of $\Omega^3_7$ given a choice of basis of $(1,0)$ vector fields.

\subsubsection*{The case of $2$-forms}

We now consider the case of a complex valued $2$-form $\beta$ on $M^7$. We can decompose it as
\begin{equation*}
\beta = \beta_2 + dt \wedge \beta_1 
\end{equation*}
where $\beta_1$ and $\beta_2$ are a $1$-form and a $2$-form on $X^6$, respectively. Again, strictly speaking, these forms also depend on the parameter $t$.

We write $\beta_2 = \beta_{2,0} + \beta^0_{1,1} + k \omega + \beta_{0,2}$, and $\beta_1 = \beta_{1,0} + \beta_{0,1}$, decomposing them into types determined by the K\"ahler structure of $X^6$.

\begin{prop} \label{2-14-prop}
The $2$-form $\eta$ is in $\Omega^2_{14}$ if and only if $k = 0$ and the following equations are satisfied:
\begin{align*}
& \beta_{0,1} \wedge \omega^2 + i \beta_{2,0} \wedge \bOm \, = \, 0, \\
& \beta_{1,0} \wedge \omega^2 - i \beta_{0,2} \wedge \Omega \, = \, 0.
\end{align*}
\end{prop}
\begin{proof}
The space $\Omega^2_{14}$ can be characterized as the space of $2$-forms $\beta$ such that $\beta \wedge \ps = 0$. We compute $\beta \wedge \ps$:
\begin{align*}
\beta \wedge \ps \, = \, \, & \left( \beta_{2,0} + \beta^0_{1,1} + k \omega + \beta_{0,2}) + dt \wedge (\beta_{1,0} + \beta_{0,1}) \right) \wedge \left( -dt \wedge \frac{1}{2i} ( \Omega - \bOm ) - \frac{\omega^2}{2} \right) \\ = \, \, & \frac{1}{2i} dt \wedge \beta_{2,0} \wedge \bOm - \frac{1}{2i} dt \wedge \beta_{0,2} \wedge \Omega - dt \wedge \beta_{1,0} \wedge \frac{\omega^2}{2} - dt \wedge \beta_{0,1} \wedge
\frac{\omega^2}{2} - \frac{k}{2} \omega^3,
\end{align*}
where other terms are zero due to type considerations and the fact that $\beta^0_{1,1} \wedge \omega^2$ is zero. Now collecting terms of the same type: $(3,3)$, $dt \wedge (3,2)$, and $dt \wedge (2,3)$,
gives the two equations above, and $k = 0$.
\end{proof}

Note that $k=0$ implies that the K\"ahler form $\omega$ has no component in $\Omega^2_{14}$, and is hence in $\Omega^2_7$. (This can also be seen from $\omega = \frac{\partial}{\partial t} \hk \ph$.)
Note also that there are no conditions on $\beta^0_{1,1}$ so all of these forms are in $\Omega^2_{14}$.

Before we move on to $\Omega^2_7$, as before let us consider the case when $\beta$ is real.  Then we have $\beta_{q,p} = \overline{\beta_{p,q}}$. Then we have
\begin{cor} \label{2-14-real-cor}
The {\em real} $2$-forms of type $\Omega^2_{14}$ are given by:
\begin{itemize} \setlength\itemsep{-1mm}
\item the {\em real} $2$-forms of type $(1,1)_0$ (where $(1,1)_0$ are the $(1,1)$ forms orthogonal to $\omega$), which is (pointwise) $8$-dimensional.
\item the real $2$-forms $\beta_{2,0} + \overline{\beta_{2,0}} + dt \wedge ( \beta_{1,0} + \overline{\beta_{1,0}} )$, where $i \overline{\beta_{2,0}} \wedge \Omega - \beta_{1,0} \wedge \omega^2 =0$, which
is (pointwise) $6$-dimensional. This is because for any real $(1,0)$ form $\beta_{1,0} + \overline{\beta_{0,1}}$, the equation $i \overline{\beta_{2,0}} \wedge \Omega - \beta_{1,0} \wedge \omega^2 = 0$
can be solved uniquely for $\beta_{2,0}$ due to the non-degeneracy of the holomorphic volume form $\Omega$.
\end{itemize}
\end{cor}
This gives the expected result for the (pointwise) dimension of the space of {\em real} $\Omega^2_{14}$ forms to be $8 + 6 = 14$.

We now consider the $2$-forms $\beta$ of type $\Omega^2_7$. This time the map $\beta \mapsto \ps \wedge \beta$ is an isomorphism of $\Omega^2_7$ onto $\Omega^6_7$.

\begin{prop} \label{2-7-prop}
The {\em real} $2$-forms of type $\Omega^2_7$ are the following:
\begin{itemize} \setlength\itemsep{-1mm}
\item $\spn( \omega )$, which is (pointwise) $1$-dimensional.
\item the real $2$-forms $\beta_{2,0} + \overline{\beta_{2,0}} + dt \wedge ( \beta_{1,0} + \overline{\beta_{1,0}} )$, where $(\beta_{2,0}, \beta_{1,0})$ is in the orthogonal complement of the kernel of the
linear map
\begin{equation*}
M : (\beta_{2,0}, \beta_{1,0}) \mapsto i \overline{\beta_{2,0}} \wedge \Omega - \beta_{1,0} \wedge \omega^2.
\end{equation*}
This space is (pointwise) $6$-dimensional.
\end{itemize}
This gives the expected result of $1 + 6 = 7$ for the (pointwise) dimension of the space of real $\Omega^2_7$ forms.
\end{prop}
\begin{proof}
Suppose that the real $2$-form $\beta$ is in $\Omega^2_7$. Since we already know that forms of the type $dt \wedge (1,1)_0$ are in $\Omega^2_{14}$, we can write that
\begin{equation*}
\beta = \beta_{2,0} + \overline{\beta_{2,0}} + k \omega + dt \wedge ( \beta_{1,0} + \overline {\beta_{1,0}} ),
\end{equation*}
where $k$ is some real function.  Consider the linear map $M$ from the underlying real vector space of $\Omega^{2,0} \oplus \Omega^{1,0}$ to the underlying real vector space of $\Omega^{3,2}$, given by
$M(\beta_{2,0}, \beta_{1,0}) = i \overline{\beta_{2,0}} \wedge \Omega - \beta_{1,0} \wedge \omega^2$.  We have shown in Corollary~\ref{2-14-real-cor} that $M$ has a $6$-dimensional kernel, so by the rank-nullity theorem, since the domain is $12$-dimensional, the orthogonal complement to the kernel is mapped isomorphically onto the image, a $6$-dimensional real vector space. Note that the map $M$
is (up to a non-zero constant factor) wedge product with the $\G$ $4$-form $\ps$. Since the $\Omega^2_7$ forms are precisely those which are mapped isomorphically onto $\Omega^6_7$ by wedge product with $\ps$, this gives those $\Omega^2_7$ forms which are mapped to $6$-forms of the form $dt \wedge \alpha_{3,2} + dt \wedge \overline{\alpha_{3,2}}$, a $6$-dimensional space, and the remaining $\Omega^2_7$ form is $\omega$, which is sent to a multiple of $\frac{\omega^3}{6} = \vol_6$ under wedge product with $\ps$.
\end{proof}

Here is another way to describe the real $\Omega^2_7$ forms. Recall that they are given by $X \hk \ph$, where $X$ is a real vector field on $M^7$. As before we write this as
\begin{equation*}
X = a^i \frac{\partial}{\partial z^i} + \bar a^i \frac{\partial}{\partial \bar z^i} + h \frac{\partial}{\partial t}
\end{equation*}
in terms of local complex coordinates $z^i$ on $X^6$. Then it is easy to compute that
\begin{align*}
& X \hk \ph \, = \, -h \omega + a^i \beta_i + \bar a^i \overline{\beta}_i, \\ \, \, & \text{ where } \ \beta _ i = \frac{1}{2} \left( \frac{\partial}{\partial z^i} \hk \Omega \right) + dt \wedge \left( \frac{\partial}{\partial z^i} \hk \omega \right).
\end{align*}
This gives a canonical basis of $\Omega^2_7$ given a choice of basis of $(1,0)$ vector fields.

\section{First order differential operators for $\G$ structures} \label{operatorssec}

In this section, we use both the local coordinate (indices) approach for $\G$ structures, as in~\cite{K2} and~\cite{KL}, as well as coordinate-free notation. Let $(M^7, \ph, \ps, g)$ be $7$-manifold with $\G$ structure. We want to study some natural first order differential operators on $(M^7, \ph, \ps, g)$. First, as on any Riemannian manifold, we have $\grad f$, the \emph{gradient} of a function $f$, which is the vector field
\begin{equation} \label{gradeq1}
(\grad f)^k \, = \, g^{ki} \nab{i} f,
\end{equation}
where $\nab{}$ is the Levi-Civita covariant derivative induced by the metric $g$. Invariantly, we have
\begin{equation} \label{gradeq2}
\grad f \, = \, (df)^{\sharp}.
\end{equation}
Next, also as on any Riemannian manifold, we have $\dive X$, the \emph{divergence} of a vector field $X$, which is the function
\begin{equation} \label{diveq1}
\dive X \, = \, \nab{i} X^i \, = \, g^{ij} \nab{i} X_j,
\end{equation}
where $X_j = g_{jk} X^k$ is the $1$-form metric dual to $X$. Invariantly, we have
\begin{equation} \label{diveq2}
\dive X \, = \, - d^* X^{\flat} \, = \, \st d \st X^{\flat},
\end{equation}
where $d^*$ is the formal adjoint to the exterior derivative $d$. (The identity $d^* = - \st d \st$ is true for $1$-forms on an odd-dimensional manifold.)

There is another natural first order differential operator on $M^7$, determined by the $\G$ structure, which we now proceed to define. Recall that $M$ has a \emph{cross product} $\cross$ on vector fields defined by
\begin{equation} \label{crosseq1}
\langle X \cross Y, Z \rangle \, = \, \ph(X,Y,Z),
\end{equation}
where $\langle \cdot, \cdot \rangle$ is the metric $g$ induced by $\ph$. Equivalently it is given by
\begin{equation} \label{crosseq2}
(X \cross Y)^{\flat} \, = \, Y \hk X \hk \ph \, = \, \st ( X^{\flat} \wedge Y^{\flat} \wedge \ps),
\end{equation}
where ${}^\flat$ is the musical isomorphism between vector fields and $1$-forms given by $g$. The cross product satisfies
\begin{equation} \label{crosseq3}
X \cross Y \, = \, - Y \cross X, \quad \quad \langle X \cross Y, X \rangle \, = \, 0, \quad \quad {|X \cross Y|}^2 \, = \, {|X \wedge Y|}^2 = {|X|}^2{|Y|}^2 - {\langle X, Y \rangle}^2,
\end{equation}
exactly like the cross product on an oriented Riemannian $3$-manifold. It differs from the $3$-dimensional case in the formula for the iteration of the cross product:
\begin{equation} \label{crosseq4}
X \cross (Y \cross Z) \, = \, - \langle X, Y \rangle Z + \langle X, Z \rangle Y - ( X \hk Y \hk Z \hk \ps)^{\sharp}.
\end{equation}
Note that we are using the sign convention given by~\eqref{convention1pheq} and~\eqref{convention1pseq}, which differs from the choice in~\cite{K1}. However, there is a sign error in the proof of Lemma 2.4.3 of~\cite{K1}, so in that convention the final term in equation~\eqref{crosseq4} should have a plus sign.

\begin{rmk} \label{dim3rmk}
For the cross product on an oriented $3$-manifold, equations~\eqref{crosseq1},~\eqref{crosseq2},~\eqref{crosseq3}, and~\eqref{crosseq4} hold with $\ph$ replaced by the \emph{volume form} $\vol_3$, and $\ps$ replaced by $\st_3 \vol_3 = 1$.
\end{rmk}

We also note that in index notation, the cross product can be written as
\begin{equation} \label{crosseq5}
(X \cross Y)^l \, = \, X^i Y^j \ph_{ijk} g^{kl},
\end{equation}
and the relation~\eqref{crosseq4} can be expressed as
\begin{equation} \label{crosseq6}
\ph_{ijk} \ph_{abc} g^{kc} \, = \, g_{ia} g_{jb} - g_{ib} g_{ja} - \ps_{ijab}.
\end{equation}

\subsection{The curl operator} \label{curlsec}

We can use the cross product to define another first order differential operator.

\begin{defn} \label{curldefn}
We define the \emph{curl} of a vector field $X$ to be the vector field $\curl X$ given by
\begin{equation} \label{curleq1}
(\curl X)^l \, = \, ( g^{ai} \nab{a} X^j) \ph_{ijk} g^{kl} \, = \, (\nab{a} X_b) g^{ai} g^{bj} \ph_{ijk} g^{kl}.
\end{equation}
Just as on an oriented Riemannian $3$-manifold, one can think of $\curl X$ as the vector field obtained by taking the cross product of the `vector field' $\nab{} = \nab{i} \frac{\partial}{\partial x^i}$ with the vector field $X$. From~\eqref{crosseq2}, we see that invariantly we have
\begin{equation} \label{curleq2}
(\curl X)^{\flat} \, = \, \st (d X^{\flat} \wedge \ps).
\end{equation}
In other words, up to $\G$-invariant isomorphisms, the vector field $\curl X$ is the projection onto the $\Omega^2_7$ component of the $2$-form $d X^{\flat}$. Note also that when the $\G$ structure is torsion-free (that is, when $\nab{} \ph = 0$), then we can forget about the parentheses in~\eqref{curleq1} and write unambiguously that $(\curl X)^l = \nab{a} X_b g^{ai} b^{bj} \, \ph_{ijk} g^{kl}$.
\end{defn}

\begin{rmk} \label{othersrmk}
By combining the exterior derivative and the various projections onto the irreducible $\G$-representations, one can define several more natural first order differential operators on a manifold with $\G$ structure. These are discussed in detail by Bryant in~\cite{Br2}. In this paper we will only consider the curl and later below the Dirac operator.
\end{rmk}

There are several relations between the operators $\grad$, $\dive$, and $\curl$ on a manifold $M$ with a torsion-free $\G$ structure. Before we present them, we need to recall some identities that are satisfied for torsion-free $\G$ structures. Let $X_k dx^k$ be a $1$-form on $M$. The Ricci identities say that
\begin{equation} \label{ricciidentityeq}
\nab{i} \nab{j} X_k - \nab{j} \nab{i} X_k \, = \, -R_{ijkl} g^{lm} X_m,
\end{equation}
where $R_{ijkl}$ is the Riemann curvature tensor. If we contract~\eqref{ricciidentityeq} with $g^{jk}$, we obtain
\begin{equation} \label{graddiveq}
(\grad \dive X)^{\flat} \, = \, \nab{i} ( g^{jk} \nab{j} X_k ) \, = \, g^{jk} \nab{j} \nab{i} X_k - R_{ijkl} g^{jk} g^{lm} X_m \, = \, g^{jk} \nab{j} \nab{i} X_k,
\end{equation}
where we have used the fact that $R_{ijkl} g^{jk} = R_{il}$ is the Ricci tensor, which vanishes for a torsion-free $\G$ structure. The Ricci-flatness of the metric also implies (by the Weitzenb\"ock formula) that the rough Laplacian agrees with (minus) the Hodge Laplacian $\Delta_d = d d^* + d^* d$ on $1$-forms. Specifically, if $X = X^k \frac{\partial}{\partial x^k}$ is a vector field on $M$, then
\begin{equation} \label{lapseq}
g^{ij} \nab{i} \nab{j} X_k \, = \, - \Delta_d X^{\flat} \, = \, - (d d^* + d^* d) X^{\flat}.
\end{equation}
Because $\ph$ is torsion-free, the Riemann curvature tensor $R_{ijkl}$ lies in $\mathrm{Sym}^2 (\Omega^2_{14})$ (see Corollary 4.7 in~\cite{K2}). Proposition 2.6 in~\cite{K2} now says that
\begin{equation*}
R_{ijkl} g^{ia} g^{jb} \ps_{abcd} \, = \, 2 R_{cdkl}.
\end{equation*}
Contracting the above identity with $g^{kc}$ now gives
\begin{equation} \label{Riemannpsicontracteq}
R_{ijkl} g^{ia} g^{jb} g^{kc} \ps_{abcd} \, = \, 2 R_{cdkl} g^{kc} \, = \, - 2 R_{dl} \, = \, 0.
\end{equation}
We are now ready to establish the relations between $\grad$, $\dive$, and $\curl$ on a $\G$ manifold.

\begin{prop} \label{divgradcurlprop}
Let $f$ be any function and $X$ be any vector field on a manifold $M$ with a \emph{torsion-free} $\G$ structure. The following relations hold:
\begin{align} \label{curlgradeq}
\curl  (\grad f)  & = \, 0, \\ \label{divcurleq} \dive (\curl X) & = \, 0, \\ \label{curlcurleq} \curl (\curl X) & = \, (\grad (\dive X)) + (\Delta_d X^{\flat})^{\sharp}.
\end{align}
\end{prop}
\begin{proof}
To establish~\eqref{curlgradeq}, we note that equations~\eqref{curleq2} and~\eqref{gradeq2} show that
\begin{equation*}
\st (\curl (\grad f))^{\flat} \, = \, d (df) \wedge \ps,
\end{equation*}
which vanishes since $d(df) = 0$. Note that~\eqref{curlgradeq} does not require the torsion-free hypothesis. To establish~\eqref{divcurleq}, we note that equations~\eqref{diveq2} and~\eqref{curleq2} show that
\begin{equation*}
\dive (\curl X) \, = \, \st d \st ( \st (d X^{\flat} \wedge \ps)) \, = \, \st d (d X^{\flat} \wedge \ps) \, = \, 0,
\end{equation*}
using the facts that $\st^2 = 1$, $d(d X^{\flat}) = 0$, and $d \ps = 0$. Note that~\eqref{divcurleq} only requires the $\G$ structure to be coclosed ($d \ps = 0$.) Finally, to prove~\eqref{curlcurleq}, we will use local coordinates, and we will require the full torsion-free hypothesis. Using~\eqref{curleq1}, we compute:
\begin{align*}
(\curl (\curl X))_k & = \, \nab{p} (\curl X)_q \, g^{pa} g^{qb} \ph_{abk} \\ & = \, \nab{p}  \left( \nab{\alpha} X_{\beta} g^{\alpha i} g^{\beta j} \ph_{ijq} \right) g^{pa} g^{qb} \ph_{abk} \\ & = \, (\nab{p} \nab{\alpha} X_{\beta} ) g^{\alpha i} g^{\beta j} g^{pa} (\ph_{ijq} \ph_{kab} g^{qb} ) \\ & = \, (\nab{p} \nab{\alpha} X_{\beta} ) g^{\alpha i} g^{\beta j} g^{pa} ( g_{ik} g_{ja} - g_{ia} g_{jk} - \ps_{ijka} ),
\end{align*}
where we have used~\eqref{crosseq6} in the last line above. This expression now simplifies to:
\begin{align*}
(\curl (\curl X))_k & = \, g^{\beta p} (\nab{p} \nab{k} X_{\beta} ) - g^{\alpha p} (\nab{p} \nab{\alpha} X_{k} )
 - (\nab{p} \nab{\alpha} X_{\beta} ) g^{\alpha i} g^{\beta j} g^{pa} \ps_{ijka} \\ & = \, (\grad (\dive X))_k + (\Delta_d X^{\flat})_k + (\nab{p} \nab{\alpha} X_{\beta} ) g^{pa} g^{\alpha i} g^{\beta j}  \ps_{aijk},
\end{align*}
using~\eqref{graddiveq} and~\eqref{lapseq}. Thus to prove~\eqref{curlcurleq} it remains to show that the last term above is zero. By the skew-symmetry of $\ps_{aijk}$, we can write this last term as:
\begin{align*}
(\nab{p} \nab{\alpha} X_{\beta} ) g^{pa} g^{\alpha i} g^{\beta j}  \ps_{aijk} & = \, \frac{1}{2} ( \nab{p} \nab{\alpha} X_{\beta} - \nab{\alpha} \nab{p} X_{\beta} )  g^{pa} g^{\alpha i} g^{\beta j}  \ps_{aijk} \\ & = \, - \frac{1}{2} R_{p\alpha \beta m} g^{mn} X_n g^{pa} g^{\alpha i} g^{\beta j}  \ps_{aijk} \, = \, 0,
\end{align*}
using~\eqref{ricciidentityeq} and~\eqref{Riemannpsicontracteq}. This completes the proof.
\end{proof}

\begin{rmk} \label{comparermk}
The identities in~\eqref{curlgradeq},~\eqref{divcurleq}, and~\eqref{curlcurleq} are the exact analogues of similar identities for oriented Riemannian $3$-manifolds. However, it is important to remember that in the non-torsion-free case, the second and third of these identities would have correction terms involving one derivative of $X$ multiplied by the torsion $T$ of the $\G$ structure.
\end{rmk}

\subsection{The Dirac operator} \label{diracsec}

In this section, we will often implicitly use the metric $g$ to identify vector fields and $1$-forms, to minimize notational clutter. Any $7$-manifold $M$ with $\G$ structure is necessarily orientable and spin, and there is a natural identification of the spinor bundle $\mathbb S$ (which is a rank $8$ real vector bundle on $M$) with the bundle $\R \oplus TM$ over $M$ whose sections are smooth functions on $M$ plus smooth vector fields on $M$, which we will now explain. At a point $p$ in $M$, the fibre of $\R \oplus TM$ is $\R \oplus T_p M \cong \R \oplus \R^7$, which we identify with the octonions $\Oc$ by identifying $\R = \real(\Oc)$ and $\R^7 = \imag (\Oc)$. Now for the fibre $\R \oplus \R^7 \cong \Oc$ to be a spinor space $\mathbb S$, there should exist a Clifford multiplication $\cdot$ of the $1$-forms on $\mathbb S$ satisfying the fundamental identity
\begin{equation} \label{cliffordeq}
X \cdot (Y \cdot s) + Y \cdot (X \cdot s) \, = \, - 2 \langle X, Y \rangle s,
\end{equation}
where $X$ and $Y$ are $1$-forms, $s = (f, Z)$ is a spinor (a pair consisting of a function $f$ and a vector field $Z$), and $\langle \cdot , \cdot \rangle$ is the inner product on $1$-forms induced by the Riemannian metric $g$ on $M$ coming from $\ph$.

\begin{lemma} \label{cliffordlemma}
Octonion multiplication by imaginary octonions $(\R^7)$ on full octonions $(\R^8)$ is a Clifford multiplication of $\Gamma (T^* M)$ on $C^{\infty}(M) \oplus \Gamma (TM)$. That is, the identity~\eqref{cliffordeq} is satisfied.
\end{lemma}
\begin{proof}
If $(f_k, X_k) \in \R \oplus \R^7$ are two `spinors' for $k=1,2$, then one can check that the octonion product of the two is:
\begin{equation} \label{octmulteq}
(f_1 , X_1) (f_2, X_2) \, = \, (f_1 f_2 - \langle X_1, X_2 \rangle \, , \, f_1 X_2 + f_2 X_1 + X_1 \cross X_2 ),
\end{equation}
where the cross product of two $1$-forms (equivalently vector fields using the metric) is given by Definition~\ref{curldefn}. Now let $s = (f, Z)$ in $\R^8$ be a spinor, and let $Y$ in $\R^7$ be a $1$-form. Then, defining the \emph{Clifford product} $\cdot$ to be given by octonion multiplication, we see by equation~\eqref{octmulteq} that
\begin{equation} \label{octmulteq2}
Y \cdot (f, Z) \, = \, ( - \langle Y, Z \rangle \, , \,  f Y + Y \cross Z).
\end{equation}
Composing with another Clifford multiplication by a $1$-form $X$ in $\R^7$, we get
\begin{align*}
& X \cdot (Y \cdot (f, Z)) \, = \, ( - \langle X, f Y + Y \cross Z \rangle \, , \, - \langle Y, Z \rangle X + X \cross (f Y + Y \cross Z)) \\ & = \, ( - f \langle X, Y \rangle - \ph(X, Y, Z) \, , \,  - \langle Y, Z \rangle X + f X \cross Y - \langle X, Y \rangle Z + \langle X , Z \rangle Y + \ps(X,Y,Z, \cdot) ),
\end{align*}
where we have used equations~\eqref{crosseq1} and~\eqref{crosseq4}. Now interchanging $X$ and $Y$ above and summing, the skew-symmetry of $\ph$, $\ps$, and $\cross$ gives
\begin{equation*}
X \cdot (Y \cdot (f, Z))  + Y \cdot (X \cdot (f, Z)) \, = \, -2 \langle X, Y \rangle (f, Z),
\end{equation*}
as claimed.
\end{proof}

Because $\mathbb S = \R \oplus TM$ is a spinor bundle on $M$, we can now define its Dirac operator.
\begin{defn} \label{diracdefn}
The \emph{Dirac operator} $\dirac$ is a first order differential operator from $\mathbb S$ to $\mathbb S$ defined as follows. Let $s = (f, X)$ be a section of $\mathbb S$. Then
\begin{equation} \label{diraceq}
\dirac (f, X) \, = \, dx^k \cdot (\nab{k} (f, X)) \, = \, dx^k \cdot (\nab{k}f, \nab{k}X).
\end{equation}
That is, $\dirac s$ should be thought of as Clifford multiplication on $s$ by the `$1$-form' $\nab{} = dx^k \nab{k}$.
\end{defn}
By equation~\eqref{octmulteq2}, we see that
\begin{equation} \label{diraceq2}
\begin{aligned}
\dirac (f, X) & = \, ( - \langle dx^k, (\nab{k} X_j) dx^j \rangle \, , \, (\nab{k}f) (dx^k)^{\sharp} + (dx^k \nab{k})^{\sharp} \cross X ) \\ & = \, ( - \dive X \, , \, \grad f + \curl X ),
\end{aligned}
\end{equation}
which expresses the Dirac operator in terms of $\dive$, $\grad$, and $\curl$.

\begin{rmk} \label{selfadjointrmk}
It is easy to check using~\eqref{octmulteq2} that Clifford multiplication is skew-adjoint with respect to the inner product on $\R \oplus TM$ given by $g$. That is,
\begin{equation*}
\langle X \cdot s_1, s_2 \rangle \, = \, - \langle s_1, X \cdot s_2 \rangle.
\end{equation*}
From this it follows that the Dirac operator is \emph{formally self-adjoint}: $\dirac^* = \dirac$. That is, the difference $\langle \dirac s_1, s_2 \rangle - \langle s_1, \dirac s_2 \rangle$ is a divergence, and hence the integral of it over a (compact) manifold $M$ will vanish by Stokes' theorem.
\end{rmk}
We now relate the \emph{Dirac Laplacian} $\dirac^* \! \dirac = \dirac^2$ to the Hodge Laplacian $\Delta_d$.

\begin{prop} \label{twolapsprop}
The Dirac Laplacian $\dirac^2$ and the Hodge Laplacian $\Delta_d$ are equal:
\begin{equation} \label{twolapseq}
\dirac^2 (f, X) \, = \, (\Delta_d f, (\Delta_d X^{\flat})^{\sharp}).
\end{equation}
\end{prop}
\begin{proof}
Using equation~\eqref{diraceq2}, we compute directly:
\begin{align*}
\dirac^2 (f, X) & = \, \dirac ( - \dive X \, , \, \grad f + \curl X ) \\ & = \, ( - \dive ( \grad f + \curl X ) \, , \, \grad( - \dive X ) + \curl( \grad f + \curl X) ) \\ & = \, ( - \dive(\grad f) - \dive(\curl X) \, , \, - \grad(\dive X) + \curl(\grad f) + \curl(\curl X) ) \\ & = \, (- \dive (\grad f) \, , \, (\Delta_d X^{\flat})^{\sharp} ),
\end{align*}
using Proposition~\ref{divgradcurlprop}. The claim now follows from the fact that $-\dive (\grad f) = - g^{ij} \nab{i} \nab{j} f = \Delta_d f$, which holds on any Riemannian manifold.
\end{proof}

\begin{rmk} \label{spinrmk}
This result is of course exactly what we expect, since on any spin manifold, the Weitzenb\"ock formula for the Dirac Laplacian shows that it differs from (minus) the rough Laplacian by a term involving the scalar curvature. On a torsion-free $\G$ manifold, which is Ricci-flat, the scalar curvature vanishes, so the Dirac Laplacian equals (minus) the rough Laplacian. However for any Ricci flat manifold, the usual Weitzenb\"ock formula says that (minus) the rough Laplacian equals the Hodge Laplacian on $1$-forms (and always on functions.) Therefore Proposition~\ref{twolapsprop} is just an explicit verification of this fact for torsion-free $\G$ manifolds.
\end{rmk}

\end{document}